\documentclass{amsart}
\usepackage[cp1251]{inputenc}
\usepackage[english]{babel}
\usepackage{amsmath}
\usepackage{amssymb}
\usepackage{amsfonts}

\newtheorem{thm}{Theorem}
\begin{document}

\thispagestyle{empty}

 \title[Matrix Fourier transform  ]{ Matrix Fourier transform with discontinuous coefficients
 }%

\author{O. Yaremko, E.Zhuravleva}%Указываем авторов

\address{Oleg Yaremko,Ekaterina Zhuravleva
\newline\hphantom{iii}Penza State University,% Место работы
\newline\hphantom{iii}str. Lermontov, 37, % Адрес (улица, дом, строение и т.п.)
\newline\hphantom{iii} 440038, Penza, Russia}%  Адрес (почтовый индекс, город, страна)
\email{yaremki@mail.ru}% Ваш электронный адрес для переписки

\maketitle {\small

\begin{quote}
\noindent{\bf Abstract. }  The explicit construction of direct and
inverse Fourier's vector transform with discontinuous coefficients is
presented. The technique of applying Fourier's vector transform with
discontinuous coefficients for solving problems of mathematical physics.Multidimensional integral transformations with non-separated variables for problems with discontinuous coefficients are constructed in this work. The coefficient discontinuities focused on the  of parallel hyperplanes. In this work explicit formulas for the kernels in the case of ideal coupling conditions are obtained; the basic identity of the integral transform is proved; technique of integral transforms is developed
\medskip

\end{quote} }


{Penza state university, PO box 440026, Penza, Lermontov's street, 37, Russia}

\section{Introduction}

Different representations of the solutions of the equilibrium equation
through functions of tension are used when solving problems by the variable
separation method. The required problem is taken to the solution of
differential equations of a more simple structure with the help of such
representations. Each functions of tension in these equations "is not
fastened" with others, but it enters into boundary conditions together with
the others. A.F.Ulitko [7] has offered rather effective method of research
of problems of mathematical physics - a method Eigen vector-valued
functions. This method is the vector analogue of the Fourier method.

The method of integral transformations is also an analytical method of the
decision of solution of problems theory of elasticity. The method of
integral transformations we consider and develop in this article. we come to
the most simple problem in space of images with the help of the integral
transformations (Fourier, Laplace, Hankel, etc.). The finding of the formula
of direct is the main difficulty in solving of problems of this approach.
Extensive enough bibliography of works on use of this method in problems of
the theory of elasticity is resulted in J.S.Ufljand's monography [2].
Method of the vector integral transforms of Fourier is equivalent the method
Eigen vector-valued functions, however, unlike the last it can to be applied
successfully be used, applied to the solution of problems of the theory of
elasticity in a piece-wise homogeneous medium. The theory of integral
transforms of Fourier with piece-wise constant coefficients in a scalar case
was studied by Ufljand J.S. [16], [17], Najda L.S. [11], Protsenko V. S
[12], [13], Lenjuk M. P [8], [9], [10]. The vector variant of a method
adapted for the solution of problems in piece-wise homogeneous medium is
developed by the author in [2], [19]. Unknown tension in the boundary
conditions and in the internal conditions of conjugation don't commit
splitting in a considered dynamic problem, so the application of the scalar
integral transforms of Fourier with piece-wise constant coefficients does
not lead to success. Method of the vector integral transforms of Fourier
with discontinuous coefficients is used for its solution in the present
work. Conformable theoretical bases of a method are presented in item 4 for
granted. The necessary proofs by the method of contour under the scheme
developed in [2] and [19]. The closed form solution of the dynamic problem
found in the use of this method in item 4.
Integral transforms arise in a natural way through the principle of linear superposition in constructing integral representations of solutions of linear differential equations.
First note that the structure of integral transforms with
the relevant variables are determined by the type of differential equation and the kind of media in which the problem is considered. Therefore  decision of integral transforms  are the problem for
mathematical physics piecewise-homogeneous (heterogeneous) media. It is clear this method is an effective for obtaining
the exact solution of boundary-value problems for  piecewise-homogeneous structures mathematical physics.

The author together with I.I.Bavrin has proposed integral transforms with non-separate variables for solving multidimensional problems  in the work \cite{yar}.

Let $V$ from $R^{n+1}$  be the half-space
\[
V=\left\{ {\left( {y_1 ,...,y_n ,x} \right)\in R^{n+1}:x>0}
\right\},
\]
then solution of the Dirichlet's problem for the half-space is expressed by Poisson formula takes the form:
\cite{bes}

\[
u(x,y)=\Gamma \left( {\frac{n+1}{2}} \right)\pi
^{-\frac{n+1}{2}}\int\limits_{y=0} {\frac{x}{\left[ (y-\eta )^2+x^2 \right]^{\frac{n+1}{2}}}f(\eta )d\eta } .
\]
Obviously Poisson's kernel  is the form of integral Laplace transform and therefore expansion of the function $f(y)$ for the eigenfunctions of the Laplace operator $\Delta$ is obtained from the  reproduce properties of the Poisson kernel:
$$
f(y)=\mathop {\lim }\limits_{\tau \to 0} \int\limits_0^\infty {\lambda
^{\frac{n}{2}}e^{-\lambda \tau }} \left(\frac{1}{\left( {\sqrt {2\pi }
} \right)^n}\int\limits_{R^n} {\frac{J_{\frac{n-2}{2}} \left( {\lambda
\left| {y-\eta } \right|} \right)}{\left| {y-\eta }
\right|^{^{\frac{n-2}{2}}}}} f\left( \eta \right)d\eta\right)d\lambda ,
$$
here $J_{\nu}$ is Bessel's function of order $\nu$ \cite{bes}.
 We may assume that integral transforms with non- separate variables are defined as follows \cite{yar} on the basis of this expansion:\\
direct integral Fourier transform has the form

\begin{equation} F\left[ f \right]\left( {y,\lambda } \right)=\frac{1}{\left( {\sqrt {2\pi }
} \right)^n}\int\limits_{R^n} {\frac{J_{\frac{n-2}{2}} \left( {\lambda
\left| {y-\eta } \right|} \right)}{\left| {y-\eta }
\right|^{^{\frac{n-2}{2}}}}} f\left( \eta \right)d\eta \equiv \hat {f}\left(
{y,\lambda } \right),
\end{equation}
inverse Fourier integral transform has the form
\begin{equation}
F^{-1} [\hat {f}](y)=
\mathop {\lim }\limits_{\tau \to 0} \int\limits_0^\infty {\lambda
^{\frac{n}{2}}e^{-\lambda \tau }} \hat {f}(y;\lambda )d\lambda \equiv f(y).
\end{equation}

In our case the construction of multi-dimensional analogues for integral transforms (1)-(2) with discontinuous coefficients is the purpose of this research.
\section{One-dimensional integral transforms with discontinuous coefficients}
In this paper integral transforms with discontinuous coefficients are constructed in accordance with  author's work \cite{10}. Let $\varphi \left( {x,\lambda } \right)$ and $\varphi ^\ast \left(
{x,\lambda } \right)$ be eigenfunctions of primal and dual problems
Sturm-Liouville for Fourier operator on sectionally  homogeneous axis $ I_n $,
\[
I_n =\left\{ {x:\;x\in \mathop
U\limits_{j=1}^{n+1} \left( {l_{j-1} ,l_j } \right),\;\,l_0 =-\infty
,\;\,l_{n+1} =\infty ,\;\,l_j <l_{j+1} ,\;\,j=\overline {1,n} } \right\}.
\]
Let us remark that eigenfunction $\varphi \left( {x,\lambda } \right)$,
\[
\varphi \left( {x,\lambda } \right)=\sum\nolimits_{k=2}^n {\theta \left(
{x-l_{k-1} } \right)\,\theta \left( {l_k -x} \right)\,\varphi _k \left(
{x,\lambda } \right)+}
\]
\[
+\,\theta \left( {l_1 -x} \right)\,\varphi _1 \left( {x,\lambda }
\right)+\theta \left( {x-l_n } \right)\,\varphi _{n+1} \left( {x,\lambda }
\right)
\]
is the solution of separated differential equations system
\[
\left( {a_m^2 \frac{d^2}{dx^2}+{\kern 1pt}\lambda ^2} \right)\,\varphi _m
\left( {x,\lambda } \right)=0,\;\;x\in \left( {l_m ,l_{m+1} } \right);\quad
m=1,...,n+1,
\]
by the coupling conditions
\[
\left[ {\alpha _{m1}^k \frac{d}{dx}+\beta _{m1}^k } \right]\varphi _k
=\left[ {\alpha _{m2}^k \frac{d}{dx}+\beta _{m2}^k } \right]\varphi _{k+1}
,
\]
\[
x=l_k ,\;\;k=1,...,n;\;\;m=1,2,
\]
on the boundary conditions
\[
\left. {\varphi _1 } \right|_{x=-\infty } =0\,,\;\,\left. {\;\varphi _{n+1}
} \right|_{x=\infty } =0.
\]
Similarly, the eigenfunction $\varphi ^\ast \left(
{x,\lambda } \right)$,
\[
\varphi ^\ast \left( {\xi ,\lambda } \right)=\sum\nolimits_{k=2}^n {\theta
\left( {\xi -l_{k-1} } \right)\,\theta \left( {l_k -\xi } \right)\,\varphi
_k^\ast \left( {\xi ,\lambda } \right)\,+}
\]
\[
+\theta \left( {l_1 -\xi } \right)\,\varphi _1^\ast \left( {\xi ,\lambda }
\right)+\theta \left( {\xi -l_n } \right)\,\varphi _{n+1}^\ast \left( {\xi
,\lambda } \right)
\]
is the solution  of separate differential equations system
\[
\left( {a_m^2 \frac{d^2}{dx^2}+{\kern 1pt}\lambda ^2} \right)\,\varphi
_m^\ast \left( {x,\lambda } \right)=0,\;\;x\in \left( {l_m ,l_{m+1} }
\right);\quad m=1,...,n+1,
\]
by the coupling conditions
\[
\frac{1}{\Delta _{1,k} }\left[ {\alpha _{m1}^k \frac{d}{dx}+\beta _{m1}^k }
\right]\varphi _k^\ast =\frac{1}{\Delta _{2,k} }\left[ {\alpha _{m2}^k
\frac{d}{dx}+\beta _{m2}^k } \right]\varphi _{k+1}^\ast ,
\quad
x=l_k ,\;\;
\]
where
\[
\Delta _{i,k} =\det \left( {{\begin{array}{*{20}c}
 {\alpha _{1i}^k } \hfill & {\beta _{1i}^k } \hfill \\
 {\alpha _{2i}^k } \hfill & {\beta _{2i}^k } \hfill \\
\end{array} }} \right)
k=1,...,n;\;\;\quad i,m=1,2,
\]
on the boundary conditions
\[
\left. {\varphi _1 } \right|_{x=-\infty } =0\,,\;\,\left. {\;\varphi _{n+1}
} \right|_{x=\infty } =0.
\]
Further normalization
eigenfunctions is accepted by the following:
\[\varphi _{n+1} \left( {x,\lambda }
\right)=e^{ia_{n+1}^{-1} x\lambda }. \quad \varphi _{n+\mbox{1}}^\ast \left(
{x,\lambda } \right)=e^{-ia_{n+1}^{-1} x\lambda }.\]
Let direct $F_{n} $ and inverse $F_{n}^{-1} $ Fourier transforms on the Cartesian axis with $ n $
division points be defined by the rules in [10] :
\begin{equation}
F_{n} \left[ f \right]\,\left( \lambda \right)=\sum\limits_{m=0}^{n+1}
{\int\limits_{l_m-1 }^{l_{m} } \; } u_{m}^\ast \left( {\xi ,\lambda }
\right)\,f_{m} \left( \xi \right)d\xi \equiv \hat {f}\left( \lambda
\right),
\end{equation}
\begin{equation}
f_k \left( x \right)=\frac{1}{\pi i}\int\limits_0^\infty {u_k \left(
{x,\lambda } \right)\hat {f}\left( \lambda \right)\lambda d\lambda .}
\end{equation}
\section{Vector Fourier transform with discontinuous coefficients}

Let's develop the method of vector Fourier transform for the solution this
problem. Let's consider Sturm--Liouville vector theory [1] about a design
bounded on the set of non-trivial solution of separate simultaneous ordinary
differential equations with constant matrix coefficients
\begin{equation}
\label{eq4}
\left( {A_m^2 \frac{d^2}{dx^2}+\lambda ^2{\rm E}+\Gamma _m^2 } \right)y_m
=0,\;\,q_m^2 =\lambda ^2{\rm E}+\Gamma _m^2 ,\;\,m=\overline {1,n+1}
\end{equation}
on the boundary conditions.
\begin{equation}
\label{eq5}
\left. {\left( {\left( {\alpha _{11}^0 +\lambda ^2\delta _{11}^0 }
\right)\frac{d}{dx}+\left( {\beta _{11}^0 +\lambda ^2\gamma _{11}^0 }
\right)} \right)y_1 } \right|_{x=l_0 } =0,\quad \left. {\left\| {y_{n+1} }
\right\|{\kern 1pt}{\kern 1pt}} \right|_{x=\infty } \,<\,\infty
\end{equation}
and conditions of the contact in the points of conjugation of intervals
\begin{equation}
\label{eq6}
\left( {\left( {\alpha _{j1}^k +\lambda ^2\delta _{j1}^k }
\right)\frac{d}{dx}+\left( {\beta _{j1}^k +\lambda ^2\gamma _{j1}^k }
\right)} \right)y_k =\left( {\left( {\alpha _{j2}^k +\lambda ^2\delta
_{j2}^k } \right)\frac{d}{dx}+} \right.\left. {\left( {\beta _{j2}^k
+\lambda ^2\gamma _{j2}^k } \right) } \right)y_{k+1} ,
\end{equation}
$x=l_k ,\;\,k=\overline {1,n} ,\;\,j=1,2.,$ where
\[
y_m \left( {x,\lambda } \right)=\left( {{\begin{array}{*{20}c}
 {y_{1m} \left( {x,\lambda } \right)} \hfill \\
 \vdots \hfill \\
 {y_{rm} \left( {x,\lambda } \right)} \hfill \\
\end{array} }} \right),
\left\| {y_m } \right\|=\sqrt {y_{1m}^2 +...+y_{rm}^2 } ,m=\overline {1,n+1}
.
\]
Let for some $\lambda $ the considered the boundary problem has a
non-trivial solution
\[
y\left( {x,\lambda } \right)=\sum\limits_{k=1}^n {\theta \left( {x-l_{k-1} }
\right)\,\theta \left( {l_k -x} \right)\,y_k \left( {x,\lambda }
\right)\,+\,\theta \left( {x-l_n } \right)\,y_{n+1} \left( {x,\lambda }
\right)} .
\]
The number $\lambda $ is called an Eigen value in this case, and the
corresponding decision $y\left( {x,\lambda } \right)$ is called Eigen
vector-valued function.

$$\alpha _{11}^0 ,\beta _{11}^0 ,\gamma _{11}^0 ,\delta _{11}^0 ,\alpha
_{j1}^k ,\beta _{j1}^k ,\gamma _{j1}^k ,\delta _{j1}^k ,\alpha _{j2}^k
,\beta _{j2}^k ,\gamma _{j2}^k ,\delta _{j2}^k ,A_j -\;\left(
{j=1,2;\;\,m=1,n+1;\;\,k=1,n} \right)$$
are matrixes of the size $r\times r$.
We will required invertible
\begin{equation}
\label{eq7}
\det \;\;M_{mk} \ne 0,\;\;\lambda \in \left. {\left[ {0,\infty } \right.}
\right)
\end{equation}
for matrixes
\[
M_{mk} \equiv \left( {{\begin{array}{*{20}c}
 {\beta _{1m}^k +\lambda ^2\gamma _{1m}^k } \hfill & {\alpha _{1m}^k
+\lambda ^2\delta _{1m}^k } \hfill \\
 {\beta _{2m}^k +\lambda ^2\gamma _{2m}^k } \hfill & {\alpha _{2m}^k
+\lambda ^2\delta _{2m}^k } \hfill \\
\end{array} }} \right),\;\,m=1,2;\;\,k=\overline {1,n} .
\]
Matrixes $A_m^2 $ and $\Gamma _m^2 $ , are is $m=\overline {1,n+1} $
-positive-defined [6]. We denote
\[
\Phi _{n+1} \left( x \right)=e^{q_{n+1} xi};\;\,\Psi _{n+1} \left( x
\right)=e^{-q_{n+1} xi};\;\,q_{n+1}^2 =A_{n+1}^{-2} \left( {\lambda ^2{\rm
E}+\Gamma ^2} \right).
\]
Define the induction relations the others n-pairs a matrix-importance
functions $\left( {\Phi _k ,\Psi _k } \right),\;\;k=1,n:$
\[
\left[ {\left( {\alpha _{j1}^k +\lambda ^2\delta _{j1}^k }
\right)\frac{d}{dx}+\left( {\beta _{j1}^k +\lambda ^2\gamma _{j1}^k }
\right)} \right]\,\left( {\Phi _k ,\Psi _k } \right)=
\]
\begin{equation}
\label{eq8}
=\left[ {\left( {\alpha _{j2}^k +\lambda ^2\delta _{j2}^k }
\right)\frac{d}{dx}+\left( {\beta _{j2}^k +\lambda ^2\gamma _{j2}^k }
\right)} \right]\,\left( {\Phi _{k+1} ,\Psi _{k+1} } \right),\quad
k=\overline {1,n} ,\quad j=\overline {1,2} .
\end{equation}
Let us introduce the following notation
\[
\left. {\mathop \Phi \limits^0 _1 \left( \lambda \right)=\left[ {\left(
{\alpha _{11}^0 +\lambda ^2\delta _{11}^0 } \right)\frac{d}{dx}+\left(
{\beta _{11}^0 +\lambda ^2\gamma _{11}^0 } \right)} \right]\Phi _1 \left(
{x,\lambda } \right)\,} \right|_{x=l_0 } ,
\]
\[
\left. {\mathop \Psi \limits^0 _1 \left( \lambda \right)=\left[ {\left(
{\alpha _{11}^0 +\lambda ^2\delta _{11}^0 } \right)\frac{d}{dx}+\left(
{\beta _{11}^0 +\lambda ^2\gamma _{11}^0 } \right)} \right]\Psi _1 \left(
{x,\lambda } \right)\,} \right|_{x=l_0 } ,
\]
\[
\Omega _k =\left( {{\begin{array}{*{20}c}
 {\Phi _k } \hfill & {\Psi _k } \hfill \\
 {\Phi _k^/ } \hfill & {\Psi _k^/ } \hfill \\
\end{array} }} \right),\quad i=\overline {1,n+1} .
\]
\begin{thm}The spectrum of the problem (\ref{eq4}),(\ref{eq5}),(\ref{eq6}) is a continuous and fills
all semi axis $\left( {0,\infty } \right)$. Sturm--Liouville theory r time
is degenerate. To each Eigen value $\lambda $ corresponds to exactly $r$ linearly
independent vector-valued functions. As the last it is possible to take $r$
columns matrix-importance functions.
\[
u\left( {x,\lambda } \right)=\sum\limits_{k=1}^n {\theta \left( {x-l_{k-1} }
\right)\,\theta \left( {l_k -x} \right)\,u_k \left( {x,\lambda }
\right)\,+\,\theta \left( {x-l_n } \right)\,u_{n+1} \left( {x,\lambda }
\right)} ,
\]
\begin{equation}
\label{eq9}
u_j \left( {x,\lambda } \right)=\Phi _j \left( {x,\lambda } \right)\mathop
{\Phi _1^{-1} }\limits^{0\;\;\;} \left( \lambda \right)-\Psi _j \left(
{x,\lambda } \right)\mathop {\Psi _1^{-1} }\limits^{0\;\;\;} \left( \lambda
\right).
\end{equation}
That is
\[
y^m\left( {x,\lambda } \right)=\left( {{\begin{array}{*{20}c}
 {u_{1m} \left( {x,\lambda } \right)} \hfill \\
 \vdots \hfill \\
 {u_{rm} \left( {x,\lambda } \right)} \hfill \\
\end{array} }} \right).
\]
\end{thm}
Dual Sturm--Liouville theory consists in a finding of the non-trivial
solution of separate simultaneous ordinary differential equations with
constant matrix coefficients.
\begin{equation}
\label{eq10}
\left( {A_m^2 \frac{d^2}{dx^2}+\lambda ^2{\rm E}+\Gamma _m^2 } \right)y_m
=0,\;\,q_m^2 =\lambda ^2{\rm E}+\Gamma _m^2 ,\;\,m=\overline {1,n+1}
\end{equation}
on the boundary conditions
\begin{equation}
\label{eq11}
\left. {\left( {\frac{d}{dx}y_1^\ast \left( {\beta _{11}^0 +\lambda ^2\gamma
_{11}^0 } \right)^{-1}+y_1^\ast \left( {\alpha _{11}^0 +\lambda ^2\delta
_{11}^0 } \right)^{-1}} \right){\kern 1pt}} \right|_{x=l_0 } =0,\quad
\;\;\left\| {y_{n+1}^\ast } \right\|\,<\,\infty ,
\end{equation}
and conditions of the contact in the points of conjugation of intervals
\[
\left( {-\frac{d}{dx}y_k^\ast ,y_k^\ast } \right)\left(
{{\begin{array}{*{20}c}
 {\beta _{11}^k +\lambda ^2\gamma _{11}^k } \hfill & {\alpha _{11}^k
+\lambda ^2\delta _{11}^k } \hfill \\
 {\beta _{21}^k +\lambda ^2\gamma _{21}^k } \hfill & {\alpha _{21}^k
+\lambda ^2\delta _{21}^k } \hfill \\
\end{array} }} \right)^{-1}=
\]
\begin{equation}
\label{eq12}
=\left( {-\frac{d}{dx}y_{k+1}^\ast ,y_{k+1}^\ast } \right)\left(
{{\begin{array}{*{20}c}
 {\beta _{12}^k +\lambda ^2\gamma _{12}^k } \hfill & {\alpha _{12}^k
+\lambda ^2\delta _{12}^k } \hfill \\
 {\beta _{22}^k +\lambda ^2\gamma _{22}^k } \hfill & {\alpha _{22}^k
+\lambda ^2\delta _{22}^k } \hfill \\
\end{array} }} \right)^{-1},\quad \quad x=l_k ,\quad k=\overline {1,n} .
\end{equation}
The solution of the boundary value problem we write in the form of
\[
y^\ast \left( {\xi ,\lambda } \right)=\sum\limits_{k=2}^n {\theta \left(
{\xi -l_{k-1} } \right)\,\theta \left( {l_k -\xi } \right)\,y_k^\ast \left(
{\xi ,\lambda } \right)\,+\theta \left( {l_1 -\xi } \right)\,y_1^\ast \left(
{\xi ,\lambda } \right)+\theta \left( {\xi -l_n } \right)\,y_{n+1}^\ast
\left( {\xi ,\lambda } \right)} ,
\]
\[
y_m^\ast \left( {\xi ,\lambda } \right)=\left( {{\begin{array}{*{20}c}
 {y_{m1}^\ast \left( {\xi ,\lambda } \right)} \hfill & \cdots \hfill &
{y_{mr}^\ast \left( {\xi ,\lambda } \right)} \hfill \\
\end{array} }} \right),
\]
\[
\left\| {y_m^\ast } \right\|=\sqrt {\left( {y_{1m}^\ast }
\right)^2+...+\left( {y_{rm}^\ast } \right)^2} ,m=\overline {1,n+1} .
\]
\begin{thm}The spectrum of the problem (\ref{eq4}),(\ref{eq5}),(\ref{eq6}) is a continuous and fills
semi axis $\left( {0,\infty } \right)$. Sturm--Liouville theory r time is
degenerate. To each Eigen value $\lambda $ corresponds to exactly $r$ linearly
independent vector-valued functions. As the last it is possible to take $r$
rows matrix-importance functions.
\[
u^\ast \left( {x,\lambda } \right)=\sum\limits_{k=1}^n {\theta \left(
{x-l_{k-1} } \right)\,\theta \left( {l_k -x} \right)\,u_k^\ast \left(
{x,\lambda } \right)\,+\,\theta \left( {x-l_n } \right)\,u_{n+1}^\ast \left(
{x,\lambda } \right)} ,
\]
\[
u_j^\ast \left( {x,\beta } \right)=\left( {\mathop \Phi \limits^0 _1 \left(
\beta \right),\mathop \Psi \limits^0 _1 \left( \beta \right)}
\right)\,\Omega _j^{-1} \left( {x,\beta } \right)\left(
{{\begin{array}{*{20}c}
 0 \hfill \\
 {\rm E} \hfill \\
\end{array} }} \right)A_j^{-2} ,
\]
That is
\begin{equation}
\label{eq13}
y^{\ast j}\left( {\xi ,\lambda } \right)=\left( {{\begin{array}{*{20}c}
 {u_{j1}^\ast \left( {\xi ,\lambda } \right)} \hfill & \cdots \hfill &
{u_{jr}^\ast \left( {\xi ,\lambda } \right)} \hfill \\
\end{array} }} \right),
j=\overline {1,r} .
\end{equation}
\end{thm}
The existence of spectral functions $u\left( {x,\lambda } \right)$ and the
conjugate spectral function $u^\ast \left( {x,\lambda } \right)$ allows to
write the a vector decomposition theorem on the set of $I_n^+ $.

\begin{thm} Let the vector-valued function $f (x)$ is defined on $I_n^+ $
continuous, absolutely integrated and has the bounded total variation. Then
for any $x\in I_n^+ $ true formula of decomposition
\[
f\left( x \right)=-\frac{1}{\pi j}\int\limits_0^\infty {u\left( {x,\lambda }
\right)} \left(  \right.\int\limits_{l_0 }^\infty {u^\ast
\left( {\xi ,\lambda } \right)f\left( \xi \right)d\xi +}
\]
\[
+\left( {\gamma _{11}^0 f_1 \left( {l_0 } \right)+\delta _{11}^0 {f}'_1
\left( {l_0 } \right)} \right)+\sum\limits_{k=1}^n {\left( {\phi _1^0 \left(
\lambda \right),\psi _1^0 \left( \lambda \right)} \right)\,\Omega _k^{-1}
\left( {l_k ,\lambda } \right)\,M_{k1}^{-1} \left( \lambda \right)\cdot }
\]
\begin{equation}
\label{eq14}
\cdot \left. {\left\{ {\left( {{\begin{array}{*{20}c}
 {\gamma _{21}^k } \hfill & {\delta _{21}^k } \hfill \\
 {\gamma _{22}^k } \hfill & {\delta _{22}^k } \hfill \\
\end{array} }} \right)\,\left( {{\begin{array}{*{20}c}
 {f_{k+1} \left( {l_k } \right)} \hfill \\
 {{f}'_{k+1} \left( {l_k } \right)} \hfill \\
\end{array} }} \right)-\left( {{\begin{array}{*{20}c}
 {\gamma _{11}^k } \hfill & {\delta _{11}^k } \hfill \\
 {\gamma _{12}^k } \hfill & {\delta _{12}^k } \hfill \\
\end{array} }} \right)\,\left( {{\begin{array}{*{20}c}
 {f_k \left( {l_k } \right)} \hfill \\
 {{f}'_k \left( {l_k } \right)} \hfill \\
\end{array} }} \right)} \right\}} \right)\lambda d\lambda .
\end{equation}
The decomposition theorem allows to enter the direct and inverse matrix
integral Fourier transform on the real semi axis with conjugation points:
\[
F_{n+} \left[ f \right]\left( \lambda \right)=\int\limits_{l_0 }^\infty
{u^\ast \left( {\xi ,\lambda } \right)f\left( \xi \right)d\xi +}
\]
\[
+\left( {\gamma _{11}^0 f_1 \left( {l_0 } \right)+\delta _{11}^0 f_1^/
\left( {l_0 } \right)} \right)+\sum\limits_{k=1}^n {\left( {\phi _1^0 \left(
\lambda \right),\psi _1^0 \left( \lambda \right)} \right)\,\Omega _k^{-1}
\left( {l_k ,\lambda } \right)\,M_{k1}^{-1} \left( \lambda \right)\cdot }
\]
\begin{equation}
\label{eq15}
\cdot \left\{ {\left( {{\begin{array}{*{20}c}
 {\gamma _{21}^k } \hfill & {\delta _{21}^k } \hfill \\
 {\gamma _{22}^k } \hfill & {\delta _{22}^k } \hfill \\
\end{array} }} \right)\,\left( {{\begin{array}{*{20}c}
 {f_{k+1} \left( {l_k } \right)} \hfill \\
 {f_{k+1}^/ \left( {l_k } \right)} \hfill \\
\end{array} }} \right)-\left( {{\begin{array}{*{20}c}
 {\gamma _{11}^k } \hfill & {\delta _{11}^k } \hfill \\
 {\gamma _{12}^k } \hfill & {\delta _{12}^k } \hfill \\
\end{array} }} \right)\,\left( {{\begin{array}{*{20}c}
 {f_k \left( {l_k } \right)} \hfill \\
 {f_k^/ \left( {l_k } \right)} \hfill \\
\end{array} }} \right)} \right\}\equiv \tilde {f}\left( \lambda \right),
\end{equation}
\begin{equation}
\label{eq16}
F_{n+}^{-1} \left[ {\tilde {f}} \right]\,\left( x \right)=-\frac{1}{\pi
i}\int\limits_0^\infty {\lambda u\left( {x,\lambda } \right)\,\tilde
{f}\left( \lambda \right)d\lambda } \equiv f\left( x \right),
\end{equation}
when
\[
f\left( x \right)=\sum\limits_{k=1}^n {\theta \left( {l_k -x}
\right)\,\theta \left( {x-l_{k-1} } \right)\,f_k \left( x \right)\,+\theta
\left( {x-l_n } \right)\,f_{n+1} \left( x \right)} .
\]
\end{thm}
 Let's result the basic identity of
integral transform of the differential operator
\[
B=\sum\limits_{j=1}^n {\theta \left( {x-l_{j-1} } \right)\,\theta \left(
{l_j -x} \right)\left( {A_j^2 \frac{d^2}{dx^2}+\Gamma _j^2 }
\right)\,+\theta \left( {x-l_n } \right)\left( {A_{n+1}^2
\frac{d^2}{dx^2}+\Gamma _{n+1}^2 } \right)} .
\]
\begin{thm} If vector-valued function
\[
f\left( x \right)=\sum\limits_{k=1}^n {\theta \left( {x-l_{k-1} }
\right)\,\theta \left( {l_k -x} \right)f_k \left( x \right)\,+\theta \left(
{x-l_n } \right)\,f_{n+1} \left( x \right)} ,
\]
is continuously differentiated on set three times, has the limit values
together with its derivatives up to the third order inclusive
\[
f_k^{(m)} \left( {l_{k-1} } \right)=f_k^{(m)} \left( {l_{k-1} +0}
\right),\;\,m=0,1,2,3;\quad k=\overline {1,n+1}
\]
Satisfies to the boundary condition on infinity
\[
\mathop {\lim }\limits_{x\to \infty } \;\left( {u^\ast \left( {x,\lambda }
\right)\frac{d}{dx}f\left( x \right)-\frac{d}{dx}u^\ast \left( {x,\lambda }
\right)\,f\left( x \right)} \right)=0
\]
Satisfies to homogeneous conditions of conjugation (\ref{eq6}), that basic identity
of integral transform of the differential operator $B$ hold
\[
F_{n+} \left[ {B\left( f \right)} \right]\,\left( \lambda \right)=-\lambda
^2\tilde {f}\left( \lambda \right)-\left\{ {\left( {\beta _{11}^0 f_1 \left(
{l_0 } \right)+\alpha _{11}^0 f_1^/ \left( {l_0 } \right)} \right)-}
\right.
\]
\begin{equation}
\label{eq17}
\left. {-\left( {\gamma _{11}^0 A_1^2 f_1^{//} \left( {l_0 } \right)+\delta
_{11}^0 A_1^2 f_1^{///} \left( {l_0 } \right)} \right)} \right\} .
\end{equation}
\end{thm}
The proof of theorems 1,2,3,4 is spent by a method of the method of contour
integration. Similarly presented to work of the author [19].

\end{document}